\documentclass[11pt,notitlepage,twoside,a4paper]{amsart}
 \usepackage{amsfonts}

\usepackage{amsmath,amssymb,enumerate}

\usepackage{epsfig,fancyhdr,color}

\usepackage{amssymb}
\usepackage{amsmath,amsthm}
\usepackage{latexsym}
\usepackage{amscd}
\usepackage{psfrag}
\usepackage{graphicx}
\usepackage[latin1]{inputenc}
\usepackage[mathcal]{eucal}

\renewcommand{\textsc}{\textcolor{red}}

\def\scal#1#2{\langle #1, #2\rangle}
\def\R#1{\mathbb{R}^{#1}}

\DeclareMathOperator{\trace}{\mathrm{trace}}

\newtheorem{theorem}{\rm\bf Theorem}[section]
\newtheorem{proposition}[theorem]{\rm\bf Proposition}

\newtheorem*{theorem 1}{\rm\bf Proposition 1}
\newtheorem*{theorem 2}{\rm\bf Proposition 2}

\theoremstyle{definition}

\theoremstyle{remark}
\newtheorem{remark}[theorem]{\rm\bf Remark}

\def\interieur#1{\mathord{\mathop{\kern 0pt #1}\limits^\circ}}

\begin{document}

\title{A generalization of Cartan's theorem on isoparametric cubics}

\author{Vladimir G. Tkachev}
\address{Department of Mathematics, Royal Institute of Technology, Sweden}
\email{tkatchev@kth.se}

\subjclass[2000]{
Primary
53C42, 
35F20, 
17A35; 
Secondary
15A63,  
17A75  
}



\keywords{Cartan's theorem, Division algebras, Composition formulas, Quadratic maps, Eiconal equation, Minimal cubic cones}

\begin{abstract}
We generalize the well-known result of \'{E}. Cartan on isoparametric cubics by showing that  a homogeneous cubic polynomial solution of the eiconal equation  $|\nabla f|^2=9|x|^4$ must be rotationally equivalent to either $x_n^3-3x_n(x_1^2+\ldots+x_{n-1}^2)$, or to one of four  exceptional Cartan cubic polynomials in dimensions $n=5,8,14,26$.
\end{abstract}

\maketitle

\section{Introduction}

In his paper \cite{Cartan}  \'{E}. Cartan found all cubic homogeneous polynomials in $\R{n}$, $n\ge 3$, satisfying the isoparametric equations
\begin{equation}\label{isop1}
|\nabla f(x)|^2=9|x|^4,
\end{equation}
\begin{equation}\label{isop2}
\Delta f(x)=0.
\end{equation}
Amazingly, the cubic solutions of (\ref{isop1})-(\ref{isop2}) can be described by means of four real division algebras $\mathbb{F}_d$ of dimension $d$, where $\mathbb{F}_1=\mathbb{R}$ (reals), $\mathbb{F}_2=\mathbb{C}$ (complexes), $\mathbb{F}_4=\mathbb{H}$ (quaternions) and $\mathbb{F}_8=\mathbb{O}$ (octonians). Cartan proved that for a cubic solution to exist the dimension $n$ must be $5,8,14$ or $26$, i.e.
$$
n=3d+2, \quad d=1,2,4,8,
$$
and, in this case, the solution is congruent (i.e. rotationally equivalent) to one of the following polynomials (cf. \cite[p.~34]{Cartan}):
\begin{equation*}\label{CartanFormula}
\begin{split}
f_d(x)& =x_{n}^3-3x_{n}x_{n-1}^2+\frac{3}{2}x_n(X_0\bar X_0+X_1\bar X_1-2X_2\bar X_2)\\
&+\frac{3\sqrt{3}}{2}x_{n-1}(X_0\bar X_0-X_1\bar X_1)
+\frac{3\sqrt{3}}{2}((X_0X_1)X_2+\bar X_2(\bar X_1\bar X_0)),
\end{split}
\end{equation*}
where $x=(X_0,X_1,X_2,x_{n-1},x_{n})$ and vector $X_k=(x_{kd+1},\ldots,x_{kd+d})$ are identified with the corresponding elements of $\mathbb{F}_d$, $k=0,1,2$, and $\bar X$ denotes the conjugate  of $X$ in $\mathbb{F}_d$.  It is not hard to prove also that all the Cartan  polynomials are irreducible.

In dimension $n=2$ there is also a \textit{reducible} polynomial satisfying (\ref{isop1})-(\ref{isop2}),
$$
f_0(x)=x_2^3-3x_2x_1^2=\mathrm{Re} (x_2+x_1 \sqrt{-1})^3,
$$
which, though having no evident relation to the division algebras, can be thought of (at least formally) as the member of the above family corresponding to $d=0$, where all $X_i$ are supposed to be zero. It is easy to see that in higher dimensions $n\ge 3$, this polynomial $f_0(x)$  gives rise to a new family of (reducible) solutions of the eiconal equation (\ref{isop1}) \textit{alone}, namely
\begin{equation}\label{Simple}
f_0(x)=x_n^3-3x_n(x_1^2+\ldots+x_{n-1}^2).
\end{equation}
Note also that $f_0$ is not congruent to any of $f_d$ in the corresponding dimensions because all  $f_d$ are harmonic, while $\Delta f_0=6(2-n)x_n\ne 0$ for $n\ge 3$.

Our main result is the following characterization of cubic solutions of equation (\ref{isop1}) alone.

\begin{theorem}\label{th1}
Any homogeneous cubic polynomial satisfying the eiconal equation $|\nabla f|^2=9|x|^4$ is rotationally equivalent to either $x_n^3-3x_n(x_1^2+\ldots+x_{n-1}^2)$, or to one of the  exceptional Cartan cubic polynomials $f_d(x)$ in dimensions $n=5,8,14,26$. In particular, irreducible cubic solutions of (\ref{isop1}) can exist only in dimensions $n=5,8,14$ and $26$.
\end{theorem}

\begin{remark}
It is well known that for $d\ne 0$ the focal varieties $f_d=0$  are minimal cones in $\R{3d+2}$ (i.e. immersed submanifolds having zero mean curvature). So far, these four Cartan cones are the only known examples of minimal cubics besides the cubic $2x_1x_2x_3+(x_1^2-x_2^2)x_4$ (a member of Lawson's family of algebraic minimal surfaces in $\R{4}$ given in \cite{Lawson}) and two additional cubics, each in dimensions $9$ and $15$, found by Wu-yi Hsiang in \cite{Hsiang67}.  In a forthcoming paper \cite{T2009} we provide a classification of minimal cubics in $\R{n}$ and  Theorem~\ref{th1} above plays a crucial role in constructing of the so-called exceptional family of minimal cubics in $\R{3k}$.
\end{remark}

I would like to thank the referee for offering useful comments and suggestions.

\section{Symmetric composition formulas}

Recall that a composition formula of size $[r,s,m]$ over the field of real numbers (see \cite{Shapiro}) is an identity
\begin{equation*}\label{comp1}
\sum_{k=1}^m b_k^2(x,y)=|x|^2|y|^2, \quad x\in \R{r},y\in \R{s},
\end{equation*}
where $b_k(x,y)$ are real bilinear forms and $|x|^2=\scal{x}{x}$ is the usual Euclidean norm of $x$. It is well known that the existence of a composition formula of size $[r,s,m]$ is equivalent to solvability of the Hurwitz matrix equations
\begin{equation}\label{HME1}
A_i^\mathrm{t} A_i=1_s,\quad 1\le i\le r,
\end{equation}
\begin{equation}\label{HME2}
A_i^\mathrm{t} A_j+A_j^\mathrm{t} A_i=0, \quad i\ne j
\end{equation}
(see, for instance \cite{Shapiro}). Here $A_i\in \R{m\times s}$ is a matrix of size $m\times s$ with real entries,  $A^\mathrm{t}$ denote the transpose matrix, and $1_k$ stands for the unit matrix in $\R{k\times k}$. It follows from (\ref{HME1}) that $m\geq \max\{r,s\}$.

If $\max\{r,s\}=m$, say $s=m$, then the celebrated Hurwitz-Radon theorem states that a composition formula of size $[r,m,m]$ exists (equivalently, the Hurwitz matrix system of size $[r,m,m]$ is solvable) iff
\begin{equation*}\label{RHur}
r\le \rho(m),
\end{equation*}
where the Hurwitz-Radon function $\rho(m)$ is defined for positive integers $m\geq 1$ by the formula
\begin{equation}\label{foll}
\rho(m)=8a+2^b, \qquad \text{where} \;m=2^{4a+b}\cdot \mathrm{odd} , \;\; 0\leq b\le 3,
\end{equation}
and $\rho(m)=0$ otherwise. In particular, for positive integers we always have $\rho(m)\le m$ with equality only if $m=1,2,4,8$. Another useful observation is that $\rho(m)=1$ if and only if $m$ is odd.

We shall need an analogue of the Hurwitz-Radon function for symmetric solutions of (\ref{HME1})-(\ref{HME2}). Given $m\geq 1$ we define $\rho_{\mathrm{symm}}(m)$ as the maximal possible $r$ such that the Hurwitz matrix equations are solvable for symmetric matrices $A_i\in \R{m\times m}$, $i=1,\ldots,r$.
\begin{proposition}
\label{pr1}
For any $m\geq 1$,
\begin{equation}\label{rhosymm}
\rho_{\mathrm{symm}}(m)=1+\rho(\frac{m}{2}).
\end{equation}
Moreover, if $\{A_i\}_{1\le i\le r}$ is a symmetric solution of (\ref{HME1})-(\ref{HME2}) for $r=\rho_{\mathrm{symm}}(m)\geq 2$ then all the matrices are trace free:
$\trace A_i=0.$
\end{proposition}

\begin{proof}
First suppose that $\rho_{\mathrm{symm}}(m)=1$. Then $m$ must be an odd number, because otherwise $m=2k$, $k\in \mathbb{Z}$, and the following two matrices
$$
A_1=\left(
        \begin{array}{cc}
          -1_{k} & 0 \\
          0 & 1_k \\
        \end{array}
      \right),
      \qquad
A_2=\left(
        \begin{array}{cc}
         0 & 1_k \\
          1_k & 0 \\
        \end{array}
      \right),
$$
provide a symmetric solution of (\ref{HME1})-(\ref{HME2}) with $r=2$. Thus $m$ is odd and it follows from (\ref{foll}) that $\rho(m)=1$. This proves (\ref{rhosymm}) for $\rho_{\mathrm{symm}}(m)=1$.

Now let us consider the case $r:=\rho_{\mathrm{symm}}(m)\geq 2$. Then we can find a symmetric solution  $\{A_i\}_{1\le i\leq r}$ of (\ref{HME1})-(\ref{HME2}). Without loss of generality we can assume that $A_r$ has the diagonal form, say $A_r=\mathrm{diag}(a_1,\ldots,a_m)$, where $a_i\in \R{}$. Then (\ref{HME1}) implies $a_i^2=1$, that is after a suitable rotation we get
\begin{equation}\label{polar}
A_r=1_t\oplus (-1)_{m-t}, \quad 0\le t\le m.
\end{equation}

We claim that  $t(m-t)\ne0$. Indeed, if $t=0$ or $t=m$ then $A_r=\pm 1_m$, hence applying (\ref{HME2}) to  $A_i$ and $A_r$ we find $A_i^\mathrm{t}+A_i=0$, $1\le i\le r-1$, which in its turn implies $A_i=0$ because $A_i$ are symmetric. But the latter contradicts to (\ref{HME1}) for $r\ge 2$. Hence $A_r$ has eigenvalues of both signs, i.e. $1\le t\le m-1$ in (\ref{polar}).

Write the remaining $A_i$ in the block form associated with the polarization of $\R{m}$ given by (\ref{polar}),
$$
A_{i}=\left(
        \begin{array}{cc}
          C_i & E_i \\
          E_i^\mathrm{t} & D_i \\
        \end{array}
      \right), \quad 1\le i \le r-1.
$$
Here $C_i\in \R{t\times t}$ and $D_i\in \R{(m-t)\times (m-t)}$ are symmetric matrices, and $E_i\in \R{t\times (m-t)}$.  Applying again (\ref{HME2}) to $A_i$ and $A_r$ we find immediately that $C_i$ and $D_i$ are zero matrices for $i\leq r-1$. Furthermore, (\ref{HME1}) yields
\begin{equation}\label{HME3}
E_i E_i^\mathrm{t} =1_{t}, \qquad E_i^\mathrm{t} E_i=1_{m-t},
\end{equation}
and setting $1\le i,j\le r-1$ in (\ref{HME2}) we get
\begin{equation}\label{HME4}
E_i E_j^\mathrm{t}+E_j E_i^\mathrm{t} =0.
\end{equation}

Now observe that identity $E_i E_i^\mathrm{t} =1_{t}$ implies $m-t\geq t$, and similarly, $E_i^\mathrm{t} E_i =1_{m-t}$ implies $t\geq m-t$. Hence $m=2t$; in particular, $m$ is en even number. It follows that all $E_i$ are quadratic matrices and equations (\ref{HME3})-(\ref{HME4}) are equivalent to the Hurwitz matrix equations of size $[r-1,t,m-t]\equiv [r-1,\frac{m}{2},\frac{m}{2}]$. This implies by the definition of $\rho$ that $r-1\le \rho(m/2)$, i.e.
$$
\rho_{\mathrm{symm}}(m)-1\le \rho(m/2).
$$

In order to prove the inverse inequality, let us fix an even $m\geq 2$ and set
$$
r:=\rho(m/2)+1\geq 2.
$$
Let  $\{E_i\}_{1\le i\le r}$ be an arbitrary  solution of (\ref{HME1})-(\ref{HME2}) of size $[r-1,\frac{m}{2},\frac{m}{2}]$. Then it is easy to check that the symmetric matrices
\begin{equation*}\label{blockk}
A_{i}=\left(
        \begin{array}{cc}
          0 & E_i \\
          E_i^\mathrm{t} & 0 \\
        \end{array}
      \right), \quad 1\le i\le r-1,
      \qquad A_r=1_{m/2}\oplus (-1)_{m/2},
\end{equation*}
give a solution to (\ref{HME3})-(\ref{HME4}) of size $[r,m,m]$. Thus
$\rho_{\mathrm{symm}}(m)\ge r=\rho(m/2)+1$,
which finishes the proof of (\ref{rhosymm}).

The last statement of the proposition easily follows from the block form of $A_i$ and the fact that trace is invariant with respect to orthogonal transformations.

\end{proof}

%

%

\section{Proof of Theorem~\ref{th1}}

Let $f$ be any cubic polynomial satisfying (\ref{isop1}). Then $f\not \equiv 0$ and it can be brought into the normal form, i.e.
\begin{equation}\label{normal}
f(x)=x_n^3+3x_nA(\bar x)+3B(\bar x), \quad \bar x=(x_1,\ldots,x_{n-1}),
\end{equation}
where $A$ is a quadratic form and $B$ is a cubic form in $\bar x$.
Indeed, the maximum value of $f(x)$ on the unit sphere $|x|=1$ is strictly  positive and attained at some point $x^0$. Then $\nabla f(x^0)=cx^0$, hence by homogeneity of $f$,
$$
c=\scal{x^0}{\nabla f(x^0)}=3f(x^0)\ne 0,
$$
and it is easily shown that in new orthogonal coordinates with $x_0$ being the $n$th vector, $f$ takes the form (\ref{normal}).

Equating $|\nabla f|^2$ to $9|x|^2$ yields
$$
x_n^2(2A+|\nabla A|^2-2|\bar x|^2)+2x_n\scal{\nabla A}{\nabla B}+(A^2+|\nabla B|^2-|\bar x|^4)=0,
$$
where $\bar x=(x_1,\ldots,x_{n-1})$. Thus
\begin{equation}\label{new1}
2A+|\nabla A|^2-2|\bar x|^2=0,
\end{equation}
\begin{equation}\label{new2}
\scal{\nabla A}{\nabla B}=0,
\end{equation}
\begin{equation}\label{new3}
A^2+|\nabla B|^2=|\bar x|^4.
\end{equation}
We can assume without loss of generality that $A$ is given in the diagonal form, say $A(\bar x)=\mathrm{diag}(a_1,\ldots,a_{n-1})$, so that (\ref{new1}) yields $2a_i^2+a_i-1=0$. This implies that $a_i$ is either $\frac{1}{2}$ or $-1$. We re-denote the coordinates such that \begin{equation}\label{contr}
A(\bar x)=\frac{1}{2}\sum_{i=1}^p \xi_i^2-\sum_{j=1}^q \eta_i^2,\qquad p+q=n-1, \quad \bar x=(\xi,\eta).
\end{equation}
Denote by $V_1$ and $V_2$ the corresponding eigenspaces of dimensions $p$ and $q$ respectively. Thus obtained polarization  $V\equiv \R{n-1}=V_1\oplus V_2$ induces the corresponding decompositions in the tensor products, in particular,
$$
{V^{*}}^{\otimes3} \simeq  V^{3,0}\oplus V^{2,1}\oplus V^{1,2}\oplus V^{0,3}, \quad V^{i,j} ={V_1^{*}}^{\otimes i}\otimes {V_2^{*}}^{\otimes (3-i)}.
$$
According to the latter decomposition we have for the cubic form $B$
$$
B=B_{3,0}+B_{2,1}+B_{1,2}+B_{0,3},
$$
where $B_{i,3-i}\equiv B_{i,3-i}(\xi,\eta)\in V^{i,3-i}$ are linear independent cubic forms. By homogeneity one  finds
\begin{equation*}
\begin{split}
\scal{\xi}{\nabla B_{i,3-i}}&=iB_{i,3-i}, \\
\scal{\eta}{\nabla B_{i,3-i}}&=(3-i)B_{i,3-i},
\end{split}
\end{equation*}
hence by virtue of (\ref{new2}),
$$
\scal{\nabla A}{\nabla B}=\scal{\xi-2\eta}{\sum_{i=0}^3\nabla B_{i,3-i}}=-6B_{3,0}-3B_{1,2}+3B_{0,3}=0.
$$
It follows from linear independence of $B_{i,3-i}$ that $B_{3,0}=B_{1,2}=B_{0,3}=0$. Thus $B\in V^{2,1}$, i.e.
\begin{equation}\label{BQ}
B\equiv B_{2,1}=\sum_{i=1}^q \eta_i Q_i(\xi),
\end{equation}
where $Q_i(\xi)\in {V_2^{\otimes*}}^2$ is a quadratic form in $\xi$.

Note that $q\geq 1$, since otherwise we would have $B\equiv 0$ and by virtue of (\ref{new3}) $A^2=|\bar x|^4$, that would imply a contradiction to (\ref{contr}), because  $A=\frac{1}{2}|\bar \xi|^2\equiv\frac{1}{2}|\bar x|^2$ for $q=0$. Thus $q=\dim V_2\ge 1$.

Note also that if $\dim V_1=p=0$ then $B=0$ and $A=-|\eta|^2$. It is easy to check that (\ref{new1})-(\ref{new3}) turn into identities and the corresponding $f$ becomes the solution of (\ref{isop1}) in the form (\ref{Simple}).

There is only remained to treat the case when both $V_1$ and $V_2$ are non-trivial:  $\dim V_k\geq 1$, $k=1,2$. We have from (\ref{new3})
$$
\biggl(\frac{1}{2}|\xi|^2-|\eta|^2\biggr)^2+\sum_{i=1}^q Q_i^2(\xi)+|\sum_{i=1}^q \eta_i \nabla Q_i(\xi)|^2=(|\xi|^2+|\eta|^2)^2.
$$
Regarding the latter equality as an identity in $\R{}[\eta_1,\ldots,\eta_q]$, one finds
\begin{equation}\label{get1}
\sum_{i=1}^q Q_i^2(\xi)=\frac{3}{4}|\xi|^4,
\end{equation}
and
\begin{equation}\label{get2}
\scal{\nabla Q_i}{\nabla Q_j}=3\delta_{ij}|\xi|^2,
\end{equation}
where $\delta_{ij}$ is the Kronecker delta. Write $Q_i$ in  matrix form
\begin{equation}\label{QA}
Q_i(\xi)=\frac{\sqrt{2}}{3}\xi^\mathrm{t} A_i\xi,
\end{equation}
where $A_i\in \R{p\times p}$ is  symmetric. It follows then from (\ref{get2}) that the symmetric matrices $\{A_i\}_{1\le 1\le q}$ solve the Hurwitz matrix equations (\ref{HME1})-(\ref{HME2}) for $s=m=p$ and $r=q$, therefore
\begin{equation}\label{recall}
q\le \rho_{\mathrm{symm}}(p).
\end{equation}

If $q=1$ then (\ref{get1})-(\ref{get2}) immediately yields $Q_1=\frac{\sqrt{2}}{3}|\xi|^2$ (the choice of sign of $Q_1$ is immaterial because we are free to change the sign of $x_n$ in (\ref{normal})). Thus
$$
f(x)=x_n^3+\frac{3}{2}x_n(x_1^2+\ldots+x_{n-2}^2-2x_{n-1}^2)+3\sqrt{3}x_{n-1}(x_1^2+\ldots+x_{n-2}^2).
$$
But the the latter polynomial is exactly the solution (\ref{Simple}) after a suitable rotation. Namely,
$$
\begin{array}{c}
f(x_1,x_2,\ldots,x_{n-1},x_n)=f_0(x_1,x_2,\ldots,\frac{\sqrt{3}}{2}x_{n-1}+\frac{1}{2}x_n, -\frac{1}{2}x_n-\frac{\sqrt{3}}{2}x_{n}).
\end{array}
$$

Finally, let us suppose that $q\ge 2$. Then (\ref{get1}) means that
\begin{equation}\label{nonconst}
y(\zeta):=\frac{2}{\sqrt{3}}(Q_1(\zeta),\ldots, Q_q(\zeta)): \; S^{p-1}\to S^{q-1}
\end{equation}
is a quadratic map sending the unit sphere $|\zeta|=1$ to the unit sphere $|y| =1$. Note also that our assumption $q\ge 2$ implies by virtue of (\ref{get2}) that the image of $y(S^{p-1})$ in $S^{q-1}$ is distinct from a point. Then one result of P.~Yiu \cite{Yiu}  provides an obstruction for a nonconstant quadratic map to exist if the dimension $q-1$ of the target sphere is too small. More specifically, let us denote by $\sigma(k)$, $k\geq 1$, the minimal possible value of $l$ for which there exists a \textit{nonconstant} homogeneous quadratic map $S^{k}\to S^l$. Then the theorem of P.~Yiu \cite[Theorem~4]{Yiu} (see also \cite{Wood} for general polynomial maps) yields a recursive formula for  $\sigma(k)$:
$$
\sigma(2^a+b)=\left\{
            \begin{array}{ll}
              2^a, & \;0\le b< \rho(2^a)  \\
              2^a+\sigma(b), & \; \rho(2^a)\le b< 2^a  \\
            \end{array}
          \right.
$$
We shall need only two easy consequences of the Yiu formula, namely, that $\sigma(m)$ is a  non-decreasing function on $\mathbb{Z}^+$, and
\begin{equation}\label{eYiu}
\sigma(2^a)=2^a, \quad a\in \mathbb{Z}^+.
\end{equation}

In this set-up, one can rewrite the existence of a nonconstant quadratic map (\ref{nonconst}) as the lower estimate
$$
 q-1\ge \sigma(p-1).
$$
Combining this with (\ref{recall}), we get after applying Proposition~\ref{pr1} that
\begin{equation}\label{ppp}
1+\sigma(p-1)\le q\le 1+\rho(\frac{p}{2}).
\end{equation}
By our assumption $q\ge 2$, hence the right inequality in (\ref{ppp}) implies that $\rho(\frac{p}{2})\ge 1$, i.e. $p$ is even. We write this as $p=2^{\nu+1}p_0$, where $p_0$ is an odd number and $\nu\ge 0$. Then (\ref{ppp}) and the definition of $\rho$ yield that
\begin{equation}\label{ppp1}
1+\sigma(2^{\nu+1}p_0-1)\le q\le 1+\rho(2^{\nu}).
\end{equation}
Notice first that $p_0=1$, because otherwise we would have $p_0\ge 3$ and by monotonicity of $\sigma$,
$$
\sigma(2^{\nu+1}p_0-1)\ge \sigma(3\cdot 2^{\nu+1}-1)\ge \sigma(2^{\nu+2})=2^{\nu+1}.
$$
But the latter contradicts to the right inequality in (\ref{ppp1}) in view of  $\rho(2^{\nu})\le 2^\nu$. Thus $p_0=1$ and we rewrite (\ref{ppp1}) as
\begin{equation*}\label{ppp2}
1+\sigma(2^{\nu+1}-1)\le q\le 1+\rho(2^{\nu}).
\end{equation*}
Now applying $\sigma(2^{\nu+1}-1)\ge \sigma(2^{\nu})=2^{\nu}$ and $\rho(2^\nu)\le 2^\nu$, we find
$$
1+2^{\nu}\le q\le 1+2^{\nu},
$$
where the right inequality is strong for $\nu\ge 4$. Thus the only possible values are $q=2^\nu+1$, $\nu=0,1,2,3$, and an easy check shows that all the values are compatible with (\ref{ppp}). The corresponding values of $p$, $q$ and  the resulting dimension $n$ are displayed in Table~\ref{tab2}.

\begin{table}[ht]
\caption{Exceptional values of $p$ and $q$}\label{tab2}
\renewcommand\arraystretch{1.5}
\noindent\[
\begin{array}{|c|c|c|c|}
\hline
\;\nu \;& p=2^{\nu+1} & q=2^\nu+1 & n=p+q+1 \\
\hline
0 & 2 & 2 & 5 \\\hline
1 & 4 & 3 & 8 \\\hline
2 & 8 & 5 & 14 \\\hline
3 & 16 & 9 & 26 \\\hline
\end{array}
\]
\end{table}

One can
see that the last column contains exactly the dimensions of the isoparametric
cubics found by \'{E}. Cartan  that were mentioned in the Introduction. To finish the proof, it suffices only to show that for the values of $p$ and $q$ as in Table~\ref{tab2}, any  solution $f$ of (\ref{isop1}) is harmonic. But this is true because  $q=2^\nu+1$ and $p=2^{\nu+1}$, hence from (\ref{contr}) we get $\Delta A(\bar x)=p-2q=-2$, and by virtue of (\ref{normal}) and (\ref{BQ}), we find
$$
\Delta f=3(2+\Delta A(\bar x))x_n +3\Delta B(\bar x)=\sum_{i=1}^q \eta_i \Delta Q_i(\xi).
$$
On the other hand, for $q\ge 2$ the matrices $A_i$ in (\ref{QA}) by Proposition~\ref{pr1} are trace free, so that $\Delta Q_i=2\trace Q_i=0$. Thus, $\Delta f=0$ and we get (\ref{isop1})-(\ref{isop2}). Applying Cartan's theorem finishes the proof.

\bibliographystyle{amsplain}

\end{document}